\documentclass[11pt]{amsart}
\usepackage{amsmath,amsthm,epsfig,amssymb}
\title{Counting all equilateral triangles in $\{0,1,...,n\}^3$}
\author{Eugen J. Ionascu}
\curraddr{(EJI) Department of Mathematics\\ Columbus State University\\4225 University Avenue\\
Columbus, GA 31907\\
and Honorific Member of the Romanian Institute of Mathematics
``Simion Stoilow" } \email{ionascu\_eugen@colstate.edu;}
\date{\today}
\begin{document}
\def\sms{\small\scshape}
\baselineskip18pt
\newtheorem{theorem}{\hspace{\parindent}
T{\scriptsize HEOREM}}[section]
\newtheorem{proposition}[theorem]
{\hspace{\parindent }P{\scriptsize ROPOSITION}}
\newtheorem{corollary}[theorem]
{\hspace{\parindent }C{\scriptsize OROLLARY}}
\newtheorem{lemma}[theorem]
{\hspace{\parindent }L{\scriptsize EMMA}}
\newtheorem{definition}[theorem]
{\hspace{\parindent }D{\scriptsize EFINITION}}
\newtheorem{problem}[theorem]
{\hspace{\parindent }P{\scriptsize ROBLEM}}
\newtheorem{conjecture}[theorem]
{\hspace{\parindent }C{\scriptsize ONJECTURE}}
\newtheorem{example}[theorem]
{\hspace{\parindent }E{\scriptsize XAMPLE}}
\newtheorem{remark}[theorem]
{\hspace{\parindent }R{\scriptsize EMARK}}
\renewcommand{\thetheorem}{\arabic{section}.\arabic{theorem}}
\renewcommand{\theenumi}{(\roman{enumi})}
\renewcommand{\labelenumi}{\theenumi}
\newcommand{\Q}{{\mathbb Q}}
\newcommand{\Z}{{\mathbb Z}}
\newcommand{\N}{{\mathbb N}}
\newcommand{\C}{{\mathbb C}}
\newcommand{\R}{{\mathbb R}}
\newcommand{\F}{{\mathbb F}}
\newcommand{\K}{{\mathbb K}}
\newcommand{\D}{{\mathbb D}}
\def\phi{\varphi}
\def\ra{\rightarrow}
\def\sd{\bigtriangledown}
\def\ac{\mathaccent94}
\def\wi{\sim}
\def\wt{\widetilde}
\def\bb#1{{\Bbb#1}}
\def\bs{\backslash}
\def\cal{\mathcal}
\def\ca#1{{\cal#1}}
\def\Bbb#1{\bf#1}
\def\blacksquare{{\ \vrule height7pt width7pt depth0pt}}
\def\bsq{\blacksquare}
\def\proof{\hspace{\parindent}{P{\scriptsize ROOF}}}
\def\pofthe{P{\scriptsize ROOF OF}
T{\scriptsize HEOREM}\  }
\def\pofle{\hspace{\parindent}P{\scriptsize ROOF OF}
L{\scriptsize EMMA}\  }
\def\pofcor{\hspace{\parindent}P{\scriptsize ROOF OF}
C{\scriptsize ROLLARY}\  }
\def\pofpro{\hspace{\parindent}P{\scriptsize ROOF OF}
P{\scriptsize ROPOSITION}\  }
\def\n{\noindent}
\def\wh{\widehat}
\def\eproof{$\hfill\bsq$\par}
\def\ds{\displaystyle}
\def\du{\overset{\text {\bf .}}{\cup}}
\def\Du{\overset{\text {\bf .}}{\bigcup}}
\def\b{$\blacklozenge$}

\def\eqtr{{\cal E}{\cal T}(\Z) }
\def\eproofi{\bsq}

\begin{abstract} We describe a procedure of counting all
equilateral triangles in the three dimensional space whose
coordinates are allowed only in the set $\{0,1,...,n\}$. This
sequence is denoted here by $ET(n)$ and it has the entry A102698 in
``The On-Line Encyclopedia of Integer Sequences". The procedure is
implemented in Maple and its main idea is based on the results in
\cite{eji}. Using this we calculated the values $ET(n)$ for n=1..55
which are included here. Some facts and conjectures about this
sequence are stated. The main of them is that $\ds \lim_{n\to
\infty} \frac{ \ln ET(n)}{\ln n+1 }$ exists.
\end{abstract} \maketitle
\section{INTRODUCTION}
If we restrict the vertices of an equilateral triangle to be in
$\Z^3$ we obtain a typical element in $\eqtr$. It is not that hard
to see that  there are no such triangles whose vertices are
contained in the coordinate planes or any other plane parallel to
one of them. Also, the sides of a triangle in $\eqtr$ cannot be of
an arbitrary length. If one such triangle is considered a whole
family in $\eqtr$ can be generated from it that have vertices in the
same plane. Moreover, we have shown in \cite{eji} the following
theorems that we are going to use in our construction here.

\begin{theorem}\label{plane}
If the triangle $\triangle OPQ\in \eqtr$ with $O$ the origin and
$l=||\overset{\rightarrow}{OP}||$ then:\par $(i)$ the points $P$ and
$Q$ are contained in a plane of equation $ax+by+cz =0$, where
\begin{equation}\label{eq1} a^2+b^2+c^2=3d^2, \
a,b,c,d \in \Z \end{equation} and $l^2=2q$;\par $(ii)$ the side
length, $l$, is of the form $\sqrt{2(m^2-mn+n^2)}$ with $m,n\in \Z$.
\end{theorem}

It is important to be able to generate all the solutions of
(\ref{eq1}):
\begin{theorem}\label{paramofplane}  The following formulae give a three integer
parameter solution of (\ref{eq1}):
\begin{equation}\label{peqabc}
\begin{cases}
a=-x_1^2+x_2^2+x_3^2-2x_1x_2-2x_1x_3\\
b=x_1^2-x_2^2+x_3^2-2x_xx_1-2x_2x_3\\
c=x_1^2+x_2^2-x_3^2-2x_3x_1-2x_3x_2\\
d=x_1^2+x_2^2+x_3^2
\end{cases},\ \  x_1,x_2,x_3\in \Z.
\end{equation}
Moreover, every solutions of  (\ref{eq1}), $a$, $b$, $c$, $d$ can be
found with (\ref{peqabc}), simplifying if necessary by a common
divisor of $a$, $b$, $c$ and $d$ with $x_1$, $x_2$, $x_3$ no bigger
in absolute value than $\sqrt{\frac{3-\sqrt{3}}{2}}q$.
\end{theorem}

We include some general remarks about the solutions of (\ref{eq1})
which are discussed in \cite{eji}:
\begin{itemize}
  \item if we assume that $gcd(a,b,c)=1$ then all $a,b,c,d$ must be
  odd integers
  \item for every $d$ odd there exist at least one solution which is
  not trivial (i.e. $a=b=c=q$)
  \item {\bf [Gauss]} positive integer $n$ can be written as a sum of three squares iff
  $n$ is not of the form $4^k(8l+7)$ with $k,l\in \Z$ (see \cite{a}
  for an elementary proof)
\end{itemize}

Our construction depends essentially on a particular solution,
$(r,s)\in \Z^2$, of the equation:
\begin{equation}\label{rs}
2(a^2+b^2)=s^2+3r^2.
\end{equation}

It turns out that this Diophantine equation has always solutions if
 $a$, $b$, $c$ and $d$ are integers satisfying (\ref{eq1}). The
 family we have mentioned  can be described as another parametrization.

\begin{theorem}\label{generalpar} Let $a$, $b$, $c$, $d$ be odd positive
integers such that $a^2+b^2+c^2=3d^2$, with $gcd(d,c)=1$. Then a
triangle $\triangle OPQ\in \eqtr $  has its vertices in the plane of
equation $a\alpha+b\beta+c\gamma=0$ iff $P(u,v,w)$ and $Q(x,y,z)$
are given by
\begin{equation}\label{paramone}
\begin{cases}
u=m_um-n_un,\\
v=m_vm-n_vn,\\
w=m_wm-n_wn, \\
\end{cases}
\ \ \ \text{and} \ \ \
\begin{cases}
x=m_xm-n_xn,\\
y=m_ym-n_yn,\\
z=m_zm-n_zn,
\end{cases}\ \
\end{equation}
with
\begin{equation}\label{paramtwo}
\begin{array}{l}
\begin{cases}
m_x=-\frac{1}{2}[db(3r+s)+ac(r-s)]/q,\ \ \ n_x=-(rac+dbs)/q\\
m_y=\frac{1}{2}[da(3r+s)-bc(r-s)]/q,\ \ \ \ \  n_y=(das-bcr)/q\\
m_z=(r-s)/2,\ \ \ \ \ \ \ \ \ \ \ \ \ \ \ \  \ \ \ \ \  \ \ \ \ \ \ n_z=r\\
\end{cases}\\
\ \ \ \text{and}\ \ \\
\begin{cases}
m_u=-(rac+dbs)/q,\ \ \ n_u=-\frac{1}{2}[db(s-3r)+ac(r+s)]/q\\
m_v=(das-rbc)/q,\ \ \ \ \   n_v=\frac{1}{2}[da(s-3r)-bc(r+s)]/q\\
m_w=r,\ \ \ \ \ \ \ \ \ \ \ \ \ \ \ \ \ \ \  n_w=(r+s)/2\\
\end{cases}
\end{array}
\end{equation}
where $q=a^2+b^2$ and $(r,s)$ is a suitable solution of (\ref{rs})
and $m,n\in\Bbb Z$.

Moreover, the side-lengths of this triangle  are equal to
$d\sqrt{2(m^2-mn+n^2)}$.
\end{theorem}

Let us observe that one can use Theorem~\ref{generalpar} as long as
$min(gcd(q,c),gcd(q,a),gcd(d,b))=1$. We have found the following
counterexample to this property: $a=55063$, $b=2396393$,
$c=5(71)(2017)(1694953)$ and $d=3(41)(3361)(1694953)$ with
$gcd(a,d)=41$, $gcd(b,d)=3361$ and $gcd(c,d)=1694953$. However we
have calculated that this property holds true for all solutions $a$,
$b$, $c$, $d$ of (\ref{eq1}) such that $gcd(a,b,c)=1$ and with all
odd $d\le 4096$. This allows us to calculate $ET(n)$ for  $n=1,..,
64$ as we will see later.

\section{Description of the procedure and some ingredients}

The idea is based on the facts above and a few other results. One
would like to first find the side lengths of the triangles in $\eqtr
\cap \cal C_n$. This will partition these triangles into clear
classes. For this purpose we will use the
Proposition~\ref{paramofplane}. Then for a given side-length, $l$,
we need to find all the possible planes that contain triangles of
sides $l$ . This gives another criteria of sub-partition even
further these triangles. Using the parametrizations given in
Theorem~\ref{generalpar} then one finds the smallest such triangle
within a given plane and which can fit in $\cal C_n$ after a
translation. Once that is obtained we have to rotate and translate
it in all possible ways, but in a pairwise disjoint manner to fill
out $\cal C_n$. A formula for the number of all these will be given.
Finally all these numbers will be added up to make up $ET(n)$.

The first fact that we will use is the following geometric
observation.

\begin{proposition}\label{maxeqtri}
The largest side length of an equilateral whose vertices are
contained in a cube of side lengths $r$ is $r\sqrt{2}$.
\end{proposition}

\proof.\  If an equilateral triangle of side lengths $l$,  has its
vertices in the cube $[0,r]^3$, we denote by $A_1$, $A_2$, and $A_3$
the areas of the three projections of this triangle on the three
coordinate planes. It is easy to see that $A^2=A_1^2+A_2^2+A_3^2$
where $A$ is the area of the equilateral triangle given. The maximum
of the area of a triangle inscribed in a square of side lengths $r$
is easy to see that is at most $r^2/2$. Hence $A^2=3l^4/16\le
3r^4/4$. This gives $l\le r\sqrt{2}$. Certainly this happens when
the vertices are at the corners diagonally opposite on each face of
the cube.\eproof

Let us work out a concrete example example: $n=4$. Using
Proposition~\ref{maxeqtri} and the part (ii) of Theorem~\ref{plane},
the side lengths can be only $\sqrt{2}$, $\sqrt{6}$, $2\sqrt{2}$,
$\sqrt{14}$, $3\sqrt{2}$, $2\sqrt{6}$, $\sqrt{26}$ or $4\sqrt{2}$.
The $d$ values here are $1$ or $3$. Since $3(1^2)=1^2+1^2+1^2$ and
$3(3^2)=1^2+1^2+5^2$ are the only solutions of (\ref{eq1}), the
parametrizations we need in this case are, as shown in \cite{eji}:

$$\cal T_{1,1,1}=\{[(0,0,0),(m, -n,n-m), (m-n,-m,n)]:m,n\in\Z, m^2+n^2\neq 0 \}$$

and

$$\cal T_{1,1,5}=\{[(0,0,0), (4m-3n,m+3n,-m),(3m+n,-3m+4n,-n)]:m,n\in \Z, m^2+n^2\neq 0\},$$
here we used the notation $\cal T_{a,b,c}$ standing for all
triangles in $\eqtr$    having a vertex the origin and the other two
in the plane $\{(\alpha,\beta,\gamma):a\alpha+b\beta+c\gamma=0\}$.

Using the first parametrization we find the $m,n$ such that the
triangle obtained after a translation fits in $\cal C_4$ and the
side lengths are $\sqrt{2}$: $T_1=\{(1,0,0),(0,1,0),(0,0,1)\}$. This
triangle can be translated in various ways inside of $\cal C_4$, and
together with all its cube symmetries and translations contribute
with a total of $512$ in $ET(4)$. We will prove a formula that gives
the total of all these triangles generated by $T_1$ inside of $\cal
C_n$. This parametrization has to be used for all the side lengths
$\sqrt{6}$, $2\sqrt{2}$, $\sqrt{14}$, $3\sqrt{2}$, $2\sqrt{6}$,
$\sqrt{26}$  and $4\sqrt{2}$: the corresponding triangles
respectively are $T_2:=\{(1,0,2),(2,1,0),(0,2,1)\}$, $2T_1$,
$T_3:=\{(2, 0, 3), (0, 3, 2), (3, 2, 0)\}$, $3T_1$, $2T_2$,
$T_4:=\{(1,4,0),(4,0,1),(0,1,4),\}$, $4T_1$. Using the same formula
we will see the contribution of all these to $ET(4)$ is
respectively: 216, 216, 128, 64, 8, 16 and 8.

There is need to use the second parametrization too since one can
take $d=3$ to obtain the side length $3\sqrt{2}$. This gives still a
new triangle $T_5:=\{(0,0,1),(1,4,0),(4,1,0)\}$ with a total of 96
other generated by it in $\cal C_4$. Tallying all these we get
$ET(4)=1264$.

As we can see from this example, one has to derive a way of finding
how many other triangles can a given one, say $T$, generate inside
of $\cal C_n$ under all possible translations and cube symmetries
roughly speaking. We need to make this a little more precise. We are
going to assume that the given triangle, $T$, that is inside $\cal
C_{t}$ is minimal in the sense that at least one of the coordinates
of the vertices in $T$ is zero and $t$ is the smallest dimension $k$
of a cube $\cal C_k$ containing $T$.

Let us denote by $O(T)$ the orbit generated by $T$ within $\cal
C_{t}$ under all translations and cube symmetries. We also need to
introduce the standard unit vectors $e_1=(1,0,0)$, $e_2=(0,1,0)$ and
$e_3=(0,0,1)$.

It is actually surprising that in order to compute the number $f$,
of all distinct triangles generated by $T$ (union of all
translations of $O(T)$) within $\cal C_n$ in terms of the following
five variables:

\begin{enumerate}
  \item $n$ -the dimension of the cube,
  \item $t$ -the maximum of all the coordinates in $T$,
  \item $\alpha(T)$ -the cardinality of $O(T)$,
  \item $\beta(T)$ -the cardinality of $O(T)\cap [O(T)+e_1]$,
  \item $\gamma(T)$ -the cardinality of $[O(T)+e_1] \cap [O(T)+e_2] $.
\end{enumerate}

\begin{theorem}\label{calc}
The function $f(T,n)$ described above is given by
\begin{equation}\label{eq2}
f(T,n)=(n+1-t)^3\alpha-3(n+1-t)^2(n-t)\beta+3(n+1-t)(n-t)^2\gamma ,
\end{equation}
 for
all $n\ge t$.
\end{theorem}

\proof. Let us consider the cube $\cal C_s=\{0,1,...,s\}^3$ where
$s=n-t$. Clearly the number of points in this set is $(s+1)^3$. Each
point $p$ in the set $\cal C_s$ is considered in here as a vector.
So, $f=|\underset {p\in \cal  C_s}{\cup} O(T)+p|$. One essential
observation here is that $|[O(T)+p]\cap [O(T)+q]|=0$ for every $p$,
$q$ such that $||p-q||\ge 1$, where
$||p-q||=\underset{i=1,2,3}{min}(|p_i-q_i|)$. This is due to the
minimality of $T$.

Let us write the elements of $\cal C_s$ in lexicographical order:
$p_1,p_2,...,p_k$ where $k=(s+1)^3$. We look now at $C_s$ as the
three dimensional grid graph. Faces in this graph are simply unit
squares formed with vertices from $\cal C_s$.   One can look at the
cardinality of $\underset {i=1.. j}{\cup} O(T)+p_i$ and show by
induction on $j$ that this is equal to

$$j\alpha-(\# edges (\cal C_s(j))\beta +(\# faces (\cal
C_s(j))\gamma$$

\n where $\cal C_s(j)$ is the graph induced in $\cal C_s$ by the
vertices $p_1,p_2,...,p_j$.  Hence we just need to compute the
number of edges and faces in $\cal C_s$. There are eight vertices in
this graph that have degree $3$ (the corners), there are $(s-1)^3$
vertices of degree $6$, also $6(s-1)^2$ vertices with degree $5$ and
finally $12(s-1)$ of degree $4$. This gives a total of
$$\frac{1}{2}[24+6(s-1)^3+30(s-1)^2+48(s-1)]=3s(s+1)^2$$
edges. The number of faces is equal to
 $\frac{1}{2}[6s^3+6s^2]=3s^2(s+1)$.\eproof

{\bf Example:} Suppose we take $T=T_5$. Then clearly $t=4$. One can
use a symbolic calculator to find $\alpha(T)=96$, $\beta(T)=24$ and
$\gamma(T)=0$. So, the contribution of $T_5$ to a cube $\cal C_n$ if
$f(T_5,n)=96(n-3)^3-72(n-3)^2(n-4)=24n(n-3)^2.$

{\bf Remark:} These facts give a way to find lower bounds for
$ET(n)$. For instance, if we put the contribution of $T_1$ and $T_2$
we obtain $ET(n)\ge 8(2n-1)(n^2-n+1)$ for all $n\ge 2$.

To generate the side lengths we would like use a well-known result
due to Euler (see \cite{r}, pp. 568 and \cite{g}, pp. 56).

\begin{proposition}\label{ejtypet} {\bf [Euler's $6k+1$]} An integer $t$ can be written as $m^2-mn+n^2$ for some
$m,n\in \Z$ if and only if in the prime factorization of $t$, $2$
and the primes of the form $6k-1$ appear to an even exponent.
\end{proposition}


\section{The code}

Using Proposition~\ref{ejtypet} and Proposition~\ref{maxeqtri} we
have the following procedure in Maple to compute the side lengths
modulo a factor of two and the square root:

\begin{itemize}
  \item sides:=proc(n)
  \item local i,j,k,L,a,m,p,q,r,ms;
  \item $L:=\{ 1\}$;$ms:=n^2$;
  \item for i from 2 to ms do
  \item $a:=ifactors(i)$; $k:=nops(a[2])$;$r:=0$;
  \item for j from 1 to k do  $m:=a[2][j][1]$; $p:=m$ mod 6; $q:=a[2][j][2]$ mod 2;
          if  r=0 and ($m=2$ or $p=5$) and $q=1$ then $r:=1$ fi;
  \item if r=0 then $L:=L$ union $\{i\}$;fi;
  \item $od;$ $L:=convert(L,list);$  $end:$
\end{itemize}

This procedure gives for $n=10$: [1, 3, 4, 7, 9, 12, 13, 16, 19, 21,
25, 27, 28, 31, 36, 37, 39, 43, 48, 49, 52, 57, 61, 63, 64, 67, 73,
75, 76, 79, 81, 84, 91, 93, 97, 100]. This gives the corresponding
side-lengths

[$\sqrt{2}$, $\sqrt{6}$,$2\sqrt{2}$, $\sqrt{14}$, $3\sqrt{2}$,
$2\sqrt{6}$, $\sqrt{26}$, $4\sqrt{2}$, $\sqrt{38}$, $\sqrt{42}$,
$5\sqrt{2}$, $3\sqrt{6}$, $2\sqrt{14}$, $\sqrt{62}$, $6\sqrt{2}$,
$\sqrt{74}$, $\sqrt{78}$, $\sqrt{86}$, $4\sqrt{6}$, $7\sqrt{2}$,
$2\sqrt{26}$, $\sqrt{114}$, $\sqrt{122}$, $3\sqrt{14}$, $8\sqrt{2}$,
$\sqrt{134}$, $\sqrt{146}$, $5\sqrt{6}$, $2\sqrt{38}$, $\sqrt{158}$,
$9\sqrt{2}$, $2\sqrt{42}$, $\sqrt{182}$, $\sqrt{186}$, $\sqrt{194}$,
$10\sqrt{2}$]

We need a procedure that will give the odd values of $d$ that
``divide" a certain side length in the sense it is possible to write
it as $d\sqrt{m^2-mn+n^2}$ with $m,n\in \Z$:

\begin{itemize}
  \item $dkl:=proc(side)$
  \item local i,x,noft,div,y,y1,z;
  \item $x:=convert(divisors(side),list)$;$noft:=nops(x)$;$div:=\{ \}$;
  \item for i from 1 to noft do $z:=x[i]$ mod 2;
  \item if $z=1$ then $y:=side/x[i]^2$; $y1:=floor(y)$; if $y=y1$ then $div:=div$ union $\{x[i]\}$;
  fi;fi;
  \item od;
  \item convert(div,list); end:
\end{itemize}

\n For instance, if $side=\sqrt{882}$ this procedure gives
$[1,3,7,21]$ which means we have at least four possible
parametrizations that we can use to find minimal equilateral
triangles.

Next we need to find all the nontrivial solutions $[a,b,c]$ of
(\ref{eq1}), given and odd positive integer $d$, with the property
$gcd(a,b,c)=1$, $0<a\le b\le c$ which is based on an internal
procedure to solve Diophantine equation $A=X^2+Y^2$:

\begin{itemize}
  \item abcsol:=proc(q) local i,j,k,u,x,y,sol,cd; $sol:=\{ \}$;
  \item for i from 1 to q do
  \item $u:=[isolve(3d^2-i^2=x^2+y^2)]$; $k:=nops(u);$
  \item for j from 1 to k do
  \item \ \  if $rhs(u[j][1])>=i$ and $rhs(u[j][2])>=i$ then
  \item \ \ \     $cd:=gcd(gcd(i,rhs(u[j][1])),rhs(u[j][2]));$
  \item \ \ \ if $cd=1$ then $sol:=sol$ union $\{sort([i,rhs(u[j][1]),rhs(u[j][2])])\};$fi;fi;
  \item od;od; convert(sol,list); end:
\end{itemize}

\n For $d=17$, $ancsol$ finds four different solutions, [[11, 11,
25], [13, 13, 23], [1, 5, 29], [7, 17, 23]], and  in a few seconds
sends out 333 solutions for $d=2007$. One interesting solution in
this last case is
$$1937^2+1973^2+2107^2=3(2007)^2.$$

Now based on the Theorem~\ref{generalpar} we take a solution of
(\ref{eq1}) as given by the procedure above and calculate the
general parametrization:
\begin{itemize}
  \item findpar:=proc(a,b,c,m,n)
  \item local i,j,r,s,sol,mx,nx,my,ny,mu,nu,mv,nv,mz,nz,mw,nw,q,d,u,v,w,x,y,z,ef,ns,om,l;
  \item $q:=a^2+b^2;$ $sol:=convert({isolve(2q=x^2+3y^2)},list);$ $ns:=nops(sol);$ $d:=\sqrt{\frac{a^2+b^2+c^2}{3}};$
  $ef:=0;$
  \item  for i from 1 to ns do
  \item if $ef=0$ then $r:=rhs(sol[i][1]);$ $s:=rhs(sol[i][2]);$
  \item  if $s^2+3r^2=2q$ then
  $mz:=(r-s)/2;$ $nz:=r;$ $mw:=r;$ $nw:=(r+s)/2;$   $mx:=-(db(3r+s)+ac(r-s))/(2q);$
  $nx:=-(rac+dbs)/q;$ $my:=(da(3r+s)-bc(r-s))/(2q);$
  $ny:=-(rbc-das)/q;$
  $mu:=-(rac+dbs)/q;$ $nu:=-(db(s-3r)+ac(r+s))/(2q);$
  $mv:=(das-rbc)/q;$ $nv:=-(da(3r-s)+bc(r+s))/(2q);$
  \item if $mx=floor(mx)$ and $nx=floor(nx)$ and $my=floor(my)$ and
  $ny=floor(ny)$ and $mu=floor(mu)$ and $nu=floor(nu)$ and $mv=floor(mv)$ and $nv=floor(nv)$ then
  \item $u:=(mu)m-(nu)n;$ $v:=(mv)m-(nv)n;$ $w:=(mw)m-(nw)n;$
$x:=(mx)m-(nx)n;$ $y:=(my)m-(ny)n;$ $z:=(mz)m-(nz)n;$
$om:=[[u,v,w],[x,y,z]];$
  \item $ef:=1;$ fi;fi; fi; od; om; end:
\end{itemize}

For the solution, $[1,5,29]$, found earlier for the case $d=17$,
$findpar$ gives $$[[-11m-13n, -21m+20n, 4m-3n], [-24m+11n, -m+21n,
m-4n]].$$

Next, using this parametrization we would like to find if there is
any equilateral triangle in $\eqtr$ which after a translation fits
inside $\cal C_{stopp}$.

\begin{description}
  \item[1] $minimaltr$:=proc(s,a,b,c,$stopp$)
  \item[2] local i,z,u,nt,d,m,n,T,$\alpha$,$\beta$,$\gamma$,tr,out,L,tri,noft,tria,orb,avb,length,$lengthn$;
  \item[3] $d:=\sqrt{\frac{a^2+b^2+c^2}{3}};$ $z:=\frac{s}{d^2};$
  $u:=convert({isolve(z=q^2-qr+r^2)},list);$ $nt:=nops(u);$
  \item[4] for i from 1 to nt do
  \item[5] $T:=findpar(a,b,c,rhs(u[i][1]),rhs(u[i][2]));$
  \item[6] $\alpha:=min(T[1][1],T[2][1],0);$
$\beta:=min(T[1][2],T[2][2],0);$ $\gamma:=min(T[1][3],T[2][3],0);$
\item[7] $tr[i]:={[T[1][1]-\alpha,T[1][2]-\beta,T[1][3]-\gamma],
[T[2][1]-\alpha,T[2][2]-\beta,T[2][3]-\gamma],[-\alpha,-\beta,-\gamma]};$
\item[8] $\begin{array}{l}
out[i]:=max(tr[i][1][1],tr[i][1][2],tr[i][1][3],tr[i][2][1], \\
tr[i][2][2],tr[i][2][3],tr[i][3][1],tr[i][3][2],tr[i][3][3]);\end{array}$

\item[9] $od;$ $L:=sort([seq(out[i],i=1..nt)]);$  $tri:=\{ \};$
\item[10]  for i from 1 to nt do if $out[i]<= stopp$ then $tri:=tri$ union $\{tr[i]\};$ $fi;$ $
od;$
\item[11] $tri:=convert(tri,list);$  $tria:=\{  \};$
\item[12] if $nops(tri)>0$ then $noft:=nops(tri);$ $tria:=\{tri[1]\};$
$orb:=transl(tri[1]);$
\item[13] for i from 1 to $noft$ do $avb:=evalb(tri[i]\ in\ orb);$
\item[14]  if $avb=false$  then $orb:=orb$ union $transl(tri[i]);$ $tria:=tria$ union $\{tri[i]\};$
\item[15]  $fi;$ $od;$ $fi;$ $tria;$ end:
\end{description}

The minimal triangle given by this procedure for $s=17\sqrt{2}$,
$a=1$, $b=5$, $c=29$, $stopp=30$: $\{[11, 21, 0], [24, 1, 3], [0, 0,
4]\}$. The last part of the procedure is actually searching for a
set of triangles that generate all the triangles in $\eqtr$ that lay
in planes of normal $(a,b,c)$ or all other $23$ possibilities
obtained by permuting the coordinates and changing signs. The  next
procedure is used above and later in order to compute the parameters
$\alpha(T)$, $\beta(T)$ and $\gamma(T)$.

\begin{itemize}
  \item transl:=proc(T)
  \item local S,Q,i,j,k,a2,b2,c2,a,b,c,d;
  \item $Q:=convert(T,list);$ $a:=max(Q[1][1],Q[2][1],Q[3][1]);$
$b:=max(Q[1][2],Q[2][2],Q[3][2]);$
$c:=max(Q[1][3],Q[2][3],Q[3][3]);$ $d:=max(a,b,c);$ $a2:=d-a;$
$b2:=d-b;$ $c2:=d-c;$  $S:=orbit(T);$
  \item  for i from 0 to a2 do
  \item   for j from 0 to b2 do
  \item     for k from 0 to c2 do
  \item       $S:=S$ union $orbit(addvect(T,[i,j,k]));$
  \item od;od;od;S;end:
\end{itemize}

Here the procedure $addvect$ and $orbit$ are:
\begin{description}
  \item[1] addvect:=proc(T,v)  local Q,a,b,c;
  \item[2] $Q:=convert(T,list);$ $a:=v[1];$ $b:=v[2];$ $c:=v[3];$
  \item[3]
  $\{[Q[1][1]+a,Q[1][2]+b,Q[1][3]+c],[Q[2][1]+a,Q[2][2]+b,Q[2][3]+c],[Q[3][1]+a,Q[3][2]+b,Q[3][3]+c]\};$
  \item[4] end:
\end{description}

and

\begin{itemize}
  \item  $orbit:=proc(T)$
  \item  $local S,Q,T1;$
  \item $Q:=convert(T,list);$
  \item  $T1:=\{[Q[1][3],Q[1][2],Q[1][1]],[Q[2][3],Q[2][2],Q[2][1]],[Q[3][3],Q[3][2],Q[3][1]]\};$
  \item $S:=orbit1(T)$ union $orbit1(T1);$  $S;$
  \item end:
\end{itemize}

The $orbit1$ takes care of the cube symmetries:

\begin{itemize}
\item orbit1:=proc(T) local
\item i,k,T1,a,b,c,x,T2,T3,T4,T5,T6,T7,T8,T9,T10,T11,T12,T13,T14,T15,T16,T17,T18,
\item T19,T20,T21,T22,T23,T24,S,Q,d,a1,b1,c1; Q:=convert(T,list);
\item a:=max(Q[1][1],Q[2][1],Q[3][1]);a1:=min(Q[1][1],Q[2][1],Q[3][1]);
\item b:=max(Q[1][2],Q[2][2],Q[3][2]);b1:=min(Q[1][2],Q[2][2],Q[3][2]);
\item c:=max(Q[1][3],Q[2][3],Q[3][3]);c1:=min(Q[1][3],Q[2][3],Q[3][3]);
\item d:=max(a,b,c); T1:=T;
\item T2:={[Q[1][2],Q[1][3],Q[1][1]],[Q[2][2],Q[2][3],Q[2][1]],[Q[3][2],Q[3][3],Q[3][1]]};
\item T3:={[Q[1][1],Q[1][3],Q[1][2]],[Q[2][1],Q[2][3],Q[2][2]],[Q[3][1],Q[3][3],Q[3][2]]};
\item T4:={[Q[1][1],Q[1][2],d-Q[1][3]],[Q[2][1],Q[2][2],d-Q[2][3]],[Q[3][1],Q[3][2],d-Q[3][3]]};
\item T5:={[Q[1][2],Q[1][3],d-Q[1][1]],[Q[2][2],Q[2][3],d-Q[2][1]],[Q[3][2],Q[3][3],d-Q[3][1]]};
\item T6:={[Q[1][1],Q[1][3],d-Q[1][2]],[Q[2][1],Q[2][3],d-Q[2][2]],[Q[3][1],Q[3][3],d-Q[3][2]]};
\item T7:={[Q[1][1],d-Q[1][2],Q[1][3]],[Q[2][1],d-Q[2][2],Q[2][3]],[Q[3][1],d-Q[3][2],Q[3][3]]};
\item T8:={[Q[1][2],d-Q[1][3],Q[1][1]],[Q[2][2],d-Q[2][3],Q[2][1]],[Q[3][2],d-Q[3][3],Q[3][1]]};
\item T9:={[Q[1][1],d-Q[1][3],Q[1][2]],[Q[2][1],d-Q[2][3],Q[2][2]],[Q[3][1],d-Q[3][3],Q[3][2]]};
\item T10:={[d-Q[1][1],Q[1][2],Q[1][3]],[d-Q[2][1],Q[2][2],Q[2][3]],[d-Q[3][1],Q[3][2],Q[3][3]]};
\item T11:={[d-Q[1][2],Q[1][3],Q[1][1]],[d-Q[2][2],Q[2][3],Q[2][1]],[d-Q[3][2],Q[3][3],Q[3][1]]};
\item T12:={[d-Q[1][1],Q[1][3],Q[1][2]],[d-Q[2][1],Q[2][3],Q[2][2]],[d-Q[3][1],Q[3][3],Q[3][2]]};
\item T13:={[Q[1][1],d-Q[1][2],d-Q[1][3]],[Q[2][1],d-Q[2][2],d-Q[2][3]],[Q[3][1],d-Q[3][2],d-Q[3][3]]};
\item T14:={[Q[1][2],d-Q[1][3],d-Q[1][1]],[Q[2][2],d-Q[2][3],d-Q[2][1]],[Q[3][2],d-Q[3][3],d-Q[3][1]]};
\item T15:={[Q[1][1],d-Q[1][3],d-Q[1][2]],[Q[2][1],d-Q[2][3],d-Q[2][2]],[Q[3][1],d-Q[3][3],d-Q[3][2]]};
\item T16:={[d-Q[1][1],d-Q[1][2],Q[1][3]],[d-Q[2][1],d-Q[2][2],Q[2][3]],[d-Q[3][1],d-Q[3][2],Q[3][3]]};
\item T17:={[d-Q[1][2],d-Q[1][3],Q[1][1]],[d-Q[2][2],d-Q[2][3],Q[2][1]],[d-Q[3][2],d-Q[3][3],Q[3][1]]};
\item T18:={[d-Q[1][1],d-Q[1][3],Q[1][2]],[d-Q[2][1],d-Q[2][3],Q[2][2]],[d-Q[3][1],d-Q[3][3],Q[3][2]]};
\item T19:={[d-Q[1][1],Q[1][2],d-Q[1][3]],[d-Q[2][1],Q[2][2],d-Q[2][3]],[d-Q[3][1],Q[3][2],d-Q[3][3]]};
\item T20:={[d-Q[1][2],Q[1][3],d-Q[1][1]],[d-Q[2][2],Q[2][3],d-Q[2][1]],[d-Q[3][2],Q[3][3],d-Q[3][1]]};
\item T21:={[d-Q[1][1],Q[1][3],d-Q[1][2]],[d-Q[2][1],Q[2][3],d-Q[2][2]],[d-Q[3][1],Q[3][3],d-Q[3][2]]};
\item T22:={[d-Q[1][1],d-Q[1][2],d-Q[1][3]],[d-Q[2][1],d-Q[2][2],d-Q[2][3]],
[d-Q[3][1],d-Q[3][2],d-Q[3][3]]};
\item T23:={[d-Q[1][2],d-Q[1][3],d-Q[1][1]],[d-Q[2][2],d-Q[2][3],d-Q[2][1]],
[d-Q[3][2],d-Q[3][3],d-Q[3][1]]};
\item T24:={[d-Q[1][1],d-Q[1][3],d-Q[1][2]],[d-Q[2][1],d-Q[2][3],d-Q[2][2]],
[d-Q[3][1],d-Q[3][3],d-Q[3][2]]};
\item
$S:=\{T1,T2,T3,T4,T5,T6,T7,T8,T9,T10,T11,T12,T13,T14,T15,
T16,T17,T18,T19,T20,\break T21,T22,T23,T24\}$;
\item S; end:
\end{itemize}

Finally, we are ready to calculate the parameters in
Theorem~\ref{calc}. We have $\alpha(T)=transl(T)$,
$\beta(T)=inters(T)$ where

\begin{itemize}
\item inters:=proc(T) local a,b,c,Q,d,S,m,i,S1,S2; Q:=convert(T,list);
\item a:=max(Q[1][1],Q[2][1],Q[3][1]);
\item b:=max(Q[1][2],Q[2][2],Q[3][2]);
\item c:=max(Q[1][3],Q[2][3],Q[3][3]);
\item d:=max(a,b,c);S2:=transl(T);S:=convert(S2,list);m:=nops(S);S1:={};
\item for i from 1 to m do $S1:=S1$ union $\{addvect(S[i],[0,0,1])\};$ od;
\item $S2$ intersect $S1;$ end:
\end{itemize}

and $\gamma(T)=intersch(T)$ where

\begin{itemize}
\item intersch:=proc(T) local a,b,c,Q,d,S,m,i,S1,S2,S3,S4;
\item Q:=convert(T,list);
\item S2:=transl(T);S:=convert(S2,list);m:=nops(S);S1:={};
\item for i from 1 to m do
\item $S1:=S1$ union $\{addvect(S[i],[0,0,1])\};$ od; $S3:=\{\};$
\item  for i from 1 m do $S3:=S3$ union $\{addvect(S[i],[0,1,0])\};$ od;
\item nops(S1 intersect S3); end:
\end{itemize}

The Theorem~\ref{calc} is then implemented in

$f:=(n,d,\alpha,\beta,\gamma)\to
(n-d+1)^3\alpha-3(n-d+1)^2(n-d)\beta+3(n-d+1)(n-d)^2\gamma$

\begin{itemize}
\item notrincn:=proc(T,n)
\item local Q,a,b,c,x,a2,b2,c2,d,y,z,w;
\item Q:=convert(T,list);
\item a2:=max(Q[1][1],Q[2][1],Q[3][1]);
\item b2:=max(Q[1][2],Q[2][2],Q[3][2]);
\item c2:=max(Q[1][3],Q[2][3],Q[3][3]);
\item d:=max(a2,b2,c2);
\item x:=nops(transl(T));y:=nops(inters(T));w:=intersch(T);
\item $z:=f(n,d,x,y,w);$ end:
\end{itemize}

In the end one has to put together all these procedures and add the
number of triangles together.

\begin{itemize}
\item main:=proc(p,lastside,nuptols)
\item local i,j,k,s,nos,div,nod,nop,sol,x,netr,noft,l,z; netr:=nuptols;
  \item  s:=sides(p);nos:=nops(s);print(s);
 \item for i from lastside to nos do
 \item   div:=dkl(s[i]);nod:=nops(div);
\item       for j from 1 to nod do
  \item       sol:=abcsol(div[j]);nop:=nops(sol);
 \item       for k from 1 to nop do
 \item        x:=minimaltr(s[i],sol[k][1],sol[k][2],sol[k][3],p);
    \item     noft:=nops(x); if $noft>=1$ then
  \item        for l from 1 to noft do
  \item           $z:=notrincn(x[l],p);$
 \item            netr:=netr+z; print(s[i],div,sol[k],x[l],checkeq(x[l]),z,netr,i);
  \item        od; fi; od; od; od; netr; end:
\end{itemize}

The values $ET(n)$ for $n=1...55$ computed with $main$ are in given
below in increasing order: 8, 80, 368, 1264, 3448, 7792, 16176,
30696, 54216, 90104, 143576, 220328, 326680, 471232, 664648, 916344,
1241856, 1655208, 2172584, 2812664, 3598664, 4553800, 5702776,
7075264, 8705088, 10628928, 12880056, 15496616, 18523472, 22003808,
26000584, 30567400, 35756776, 41631672, 48278136, 55753272,
64134536, 73495760, 83924408, 95513248, 108379264, 122661840,
138315720, 155613408, 174622488, 195478424, 218279240, 243170376,
270288064, 299790968, 331832248, 366610560, 404253120, 444911712,
and 488902856.

\section{Some facts and conjectures}
If we look at the sequence $a_n=\frac{\ln ET(n)}{\ln (n+1)}$, $n\in
\N$ it seems like it is increasing. This sequence is clearly bounded
from above since the number of all triangles in $\{0,...,n\}^3$ is
not more than $(n+1)^9$ and so $a_n\le 9$. Numerically, the best
upper-bound for $a_n$ seems to be $5$ which is equivalent to saying
that $ET(n)\le (n+1)^5$ for all $n\in \N$.

From what we have seen before each class of triangles determined by
$a,b,c,d$, a solutions of $(\ref{eq1})$, brings in a contribution
that is a polynomial in terms of $n$. If we add these polynomials
together, we get a polynomial which can be expressed is  in the
variable $\zeta=n-1$ ($n=\zeta+1$) as follows:

\begin{description}
  \item[$n=1$]  $p_1(\zeta)=8\zeta^3+24\zeta^2+24\zeta+8$,
  $ET(1)=p_1(0)=8$;
  \item[$n=2$]
  $p_2(\zeta)=p_1(\zeta)+16+48(\zeta-1)+16(\zeta-1)^3+48(\zeta-1)^2$,
  $ET(2)=p_2(1)=80$;
  \item[$n=3$]
  $p_3(\zeta)=p_2(\zeta)+24+72(\zeta-2)+24(\zeta-2)^3+72(\zeta-2)^2$,
  $ET(3)=p_3(2)=368$;
  \item[$n=4$]
  $p_4(\zeta)=p_3(\zeta)+128+312(\zeta-3)+56(\zeta-3)^3+240(\zeta-3)^2$,
  $ET(4)=p_4(3)=1264$;
  \item[$n=5$] $p_5(\zeta)=p_4(\zeta)+40+120(
  \zeta-4)+40(\zeta-4)^3+120(\zeta-4)^2$, $ET(5)=p_5(4)=3448$;
  \item [$n=6$]
  $p_6(\zeta)=p_5(\zeta)+48+144(\zeta-5)+48(\zeta-5)^3+144(\zeta-5)^2$, $ET(6)=p_6(5)=7792$;
  \item [$n=7$] $p_7(\zeta)=p_6(\zeta)+776+1392(\zeta-6)+128(\zeta-6)^3+744(\zeta-6)^2$, $ET(7)=p_7(6)=16176$;
\item[$n=8$]
$p_8(\zeta)=p_7(\zeta)+232+552(\zeta-7)+88(\zeta-7)^3+408(\zeta-7)^2$,
$ET(8)=p_8(7)=30696$;
\item [$n=9$]
$p_9(\zeta)=p_8(\zeta)+360+840(\zeta-8)+120(\zeta-8)^3+600(\zeta-8)^2$,
$ET(9)=p_9(8)=54216$;
\item [$n=10$]
$p_{10}(\zeta)=p_9(\zeta)+80+80(\zeta-9)^3+240(\zeta-9)^2+240(\zeta-9)$,
$ET(10)=p_{10}(9)=90104$;
\end{description}
...............................................

 We conjecture that in general
\begin{equation}\label{conjecture1}
p_n(\zeta)=p_{n-1}(\zeta)+u_n(\zeta-n+1)^3+v_n(\zeta-n+1)^2+w_n(\zeta-n+1)+s_n,
n\in \N,
\end{equation}
with $u_n$,$v_n$,$w_n$, and $s_n$ non-negative integers.

\begin{center}
\epsfig{file=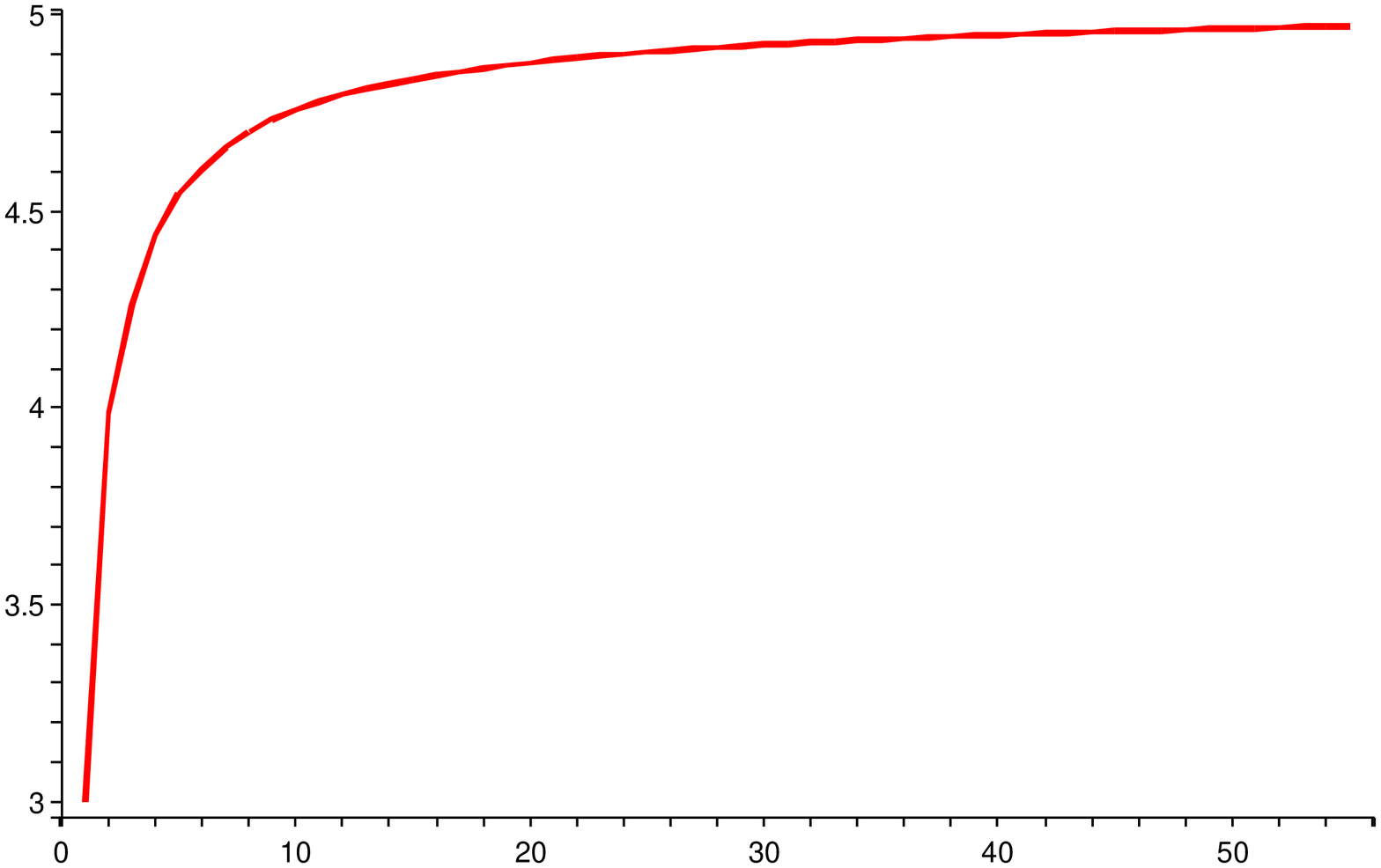,height={2in},width={3.5in}}
\end{center}

As the graph above of $n\to  \frac{\ln ET(n)}{\ln (n+1)}$ suggests,
the second conjecture is that the following limit exists

\begin{equation}\label{conjecture2}
\lim_{n\to \infty} \frac{\ln ET(n)}{\ln (n+1)} =C.
\end{equation}
  {\bf Acknowledgements:} I would like to thank F.Luca for helping
me find the counterexample mentioned after the
Theorem~\ref{generalpar}.

\end{document}